\documentclass[11pt]{article}
\usepackage{amscd}
\usepackage{amsfonts}
\usepackage{amsmath}
\usepackage[misc]{ifsym}
\usepackage{amssymb}
\usepackage{amsthm}
\usepackage{bbm}
\usepackage{CJK}
\usepackage{fancyhdr}
\usepackage{graphicx}
\usepackage{indentfirst}
\usepackage{latexsym,bm}
\usepackage{mathrsfs}
\usepackage{xypic}
\usepackage[square, comma, sort&compress, numbers]{natbib}
\usepackage[colorlinks,linkcolor=red, anchorcolor=blue]{hyperref}
\allowdisplaybreaks[4]
\newtheorem{theorem}{Theorem}[section]
\newtheorem{lemma}[theorem]{Lemma}
\newtheorem{definition}[theorem]{Definition}
\newtheorem{proposition}[theorem]{Proposition}
\newtheorem{example}[theorem]{Example}

\newtheorem{corollary}[theorem]{Corollary}

\usepackage[top=1in,bottom=1in,left=1.25in,right=1.25in]{geometry}
\def\<{\langle}
\def\>{\rangle}

\def\bu{\bullet}

\def\D{\Delta}

\def\lh{\leftharpoonup }
\def\m{\mapsto}
\def\o{\otimes}

\def\r{\rho}
\def\ra{\rightarrow}

\def\rh{\rightharpoonup }

\def\tr{\triangleright}
\def\tl{\triangleleft}

\def\v{\varepsilon}

\def\z{\zeta}

\date{}
\begin{document}
\renewcommand{\baselinestretch}{1.2}
\renewcommand{\arraystretch}{1.0}
\title{\bf Symmetries in Yetter-Drinfel'd-Long categories}
 \date{}
\author {{\bf \textbf{Dongdong Yan}$^{\dag}$, \textbf{shuanhong Wang} }\footnote{Corresponding author: E-mail: shuanhwang@seu.edu.cn}
\footnote{School of Mathematics, Southeast University, Nanjing 210096, Jiangsu, China.}}

 \maketitle
\begin{center}
\begin{minipage}{12.cm}

\noindent{\bf Abstract.} Let $H$ be a Hopf algebra and $\mathcal{LR}(H)$ the category of Yetter-Drinfel'd-Long bimodules over $H$. We first give sufficient and necessary conditions for $\mathcal{LR}(H)$ to be symmetry and pseudosymmetry, respectively. We then introduce the definition of $u$-condition in $\mathcal{LR}(H)$ and discuss the relation between the $u$-condition and the symmetry of $\mathcal{LR}(H)$. Finally, we show that $\mathcal{LR}(H)$ over a triangular (cotriangular, resp.) Hopf algebra contains a rich symmetric subcategory.
\\

\\

 \end{minipage}
 \end{center}
 \normalsize\vskip1cm

\section{Introduction}
The notion of symmetric category is a classical concept in category theory. Cohen and Westreich \cite{CW98} tested symmetries and the $u$-condition in the Yetter-Drinfel'd category $\!^{H}_{H}\mathcal{YD}$ over Hopf algebra $H$. Pareigis \cite{P01} found the necessary and sufficient condition for $\!^{H}_{H}\mathcal{YD}$ to be symmetric. Later, Panaite et al. \cite{PSV10} proposed the definition of pseudosymmetric braided categories which can be viewed as a kind of weakened symmetric braided categories, and showed that the category $\!_{H}\mathcal{YD}^{H}$ is pseudosymmetric if and only if $H$ is commutative and cocommutative. The generalization of those classical structures and results have been introduced and discussed by many authors \cite{MLC18, WG14,ZLW15}.

It is known that the Radford biproduct has a categorical interpretation (due to majid): $(H, A)$ is an admissible pair (see \cite{R85}) if and only if $A$ is a bialgebra in the Yetter-Drinfel'd category $\!^{H}_{H}\mathcal{YD}$. Panaite and Van Oystaeyen \cite{PV10} described a similar interpretation for L-R-admissible pairs and defined a prebraided category $\mathcal{LR}(H)$ (which is braided if $H$ has a bijective antipode) which contains $\!^{H}_{H}\mathcal{YD}$ and $\mathcal{YD}^{H}_{H}$ as braided subcategories. They then showed that $(H, B)$ is an L-R-admissible pair is equivalent to $D$ is a bialgebra in $\mathcal{LR}(H)$ with a extra condition
\begin{align*}
b_{(0)}\tl b'_{[-1]}\o b_{(1)}\tr b'_{[0]}=b\o b', \qquad for ~any ~b, b'\in B,
\end{align*}
where the L-R-admissible pair is the sufficient condition for L-R smash biproduct $B\bowtie H$ to be a bialgebra. The Radford biproduct is a particular case.

The aim of the present paper is to discuss the symmetries, the pseudosymmetries and the $u$-condition in Yetter-Drinfel'd-Long categories.

This paper is organized as follows: In section 1, we recall some basic definitions and results related to Yetter-Drinfel'd-Long bimodules. Then we give some examples of Yetter-Drinfel'd-Long bimodules. In section 2, we show that the Yetter-Drinfel'd-Long category $\mathcal{LR}(H)$ is symmetric if and only if $H$ is trivial in four different methods, and that $\mathcal{LR}(H)$ is pseudosymmetric if and only if $H$ is commutative and cocommutative. In section 3, we introduce the definition of $u$-condition in $\mathcal{LR}(H)$ and give a necessary and sufficient condition for $H_{i}$ $(i=1,2,3,4)$ to satisfy the $u$-condition, where $H_{i}$ is defined in Example \ref{E2.1}. Then we study the relation between the $u$-condition and the symmetry of $\mathcal{LR}(H)$. In section 4, we prove that the subcategory $\!_{H}\mathcal{M}_{H}$ of $\mathcal{LR}(H)$ over triangular Hopf algebra $H$ is symmetric. If we consider $M=H\o H$, we prove the converse. That is, assume that the braiding $\psi_{H\o H,H\o H}$ is symmetric forces $H$ to be triangular. In section 5, we give the dual cases of section 4.

\section{Preliminaries}
\def\theequation{2.\arabic{equation}}
\setcounter{equation} {0}
Throughout this paper, all algebraic systems are over a field $\Bbbk$. For a coalgebra $C$, the comultiplication will be denoted by $\D$. We follow the Sweedler's notation $\D(c)=c_{1}\o c_{2}$, for any $c\in C$, in which we often omit the summation symbols for convenience. For any vector spaces $M$ and $N$, we use $\tau: M\o N\ra N\o M$ for the flip map.

Let $A$ be a algebra, A \emph{right $A$-module} is a pair $(M, \tl)$, in which $M$ is a vector space and $\tl:M\o A\ra M$ is a linear map, called the action of $A$ on $M$, with notation $\tl (m\o a)=m\tl a$, such that, for any $a, b\in A$ and $m\in M$:
\begin{align*}
\begin{cases}
&m\tl ab=(m\tl a)\tl b,
\\&m\tl 1=m.
\end{cases}
\end{align*}

Similarly, we can define the left $A$-module. A \emph{right $A$-linear} is a linear map $f:M\ra N$ such that $f(m)\tl a=f(m\tl a)$, for any $a\in A$ and $m\in M$.

Let $C$ be a coalgebra, A \emph{right $C$-comodule} is a pair $(M, \r)$, in which $M$ is a vector space and $\r:M\ra M\o C$ is a linear map, called the coaction of $C$ on $M$, with notation $\r(m)=m_{(0)}\o m_{(1)}$, such that, for any $m\in M$:
\begin{align*}
\begin{cases}
&m_{(0)(0)}\o m_{(0)(1)}\o m_{(1)}=m_{(0)}\o m_{(1)1}\o m_{(1)2},
\\&m_{(0)}\v(m_{(1)})=m.
\end{cases}
\end{align*}

Similarly, we can define the left $C$-comodule. A \emph{right $C$-colinear} is a linear map $f:M\ra N$ such that $\r_{N}\circ f=(f\o id)\circ \r_{M}$.

Let $A$ be a algebra, and assume that $M$ are both left $A$-module via $\tr:A\o M\ra M, a\o m\m a\tr m$ and right $A$-module via $\tl: M\o A\ra M, m\o b\m m\tl b$, then $M$ is called a \emph{$A$-bimodule} if
\begin{align}\label{e2.a}
(a\tr m)\tl b=a\tr (m\tl b),
\end{align}
for any $a,b\in A$ and $m\in M$.

Let $C$ be a coalgebra, and assume that $M$ are both left $C$-comodule via $\r^{l}:M\ra C\o M, m\m m_{[-1]}\o m_{[0]}$ and right $C$-comodule via $\r^{r}:M\ra M\o C, m\m m_{(0)}\o m_{(1)}$, then $M$ is called a \emph{$C$-bicomodule} if
\begin{align}\label{e2.c}
m_{[-1]}\o m_{[0](0)}\o m_{[0](1)}=m_{(0)[-1]}\o m_{(0)[0]}\o m_{(1)},
\end{align}
for any $m\in M$.

Let $H$ be a Hopf algebra, we can denote those categories by $\!_{H}\mathcal{M}_{H}$ and $\!^{H}\mathcal{M}^{H}$. Take $\!_{H}\mathcal{M}_{H}$ whose objects are all $H$-bimodules, the morphisms in the category are morphisms of $H$-bilinear.
\begin{definition}(\cite{PV10})
Let $H$ be a Hopf algebra. A Yetter-Drinfel'd-Long bimodule over $H$ is a vector space $M$ endowed with $H$-bimodule and $H$-bicomodule structures (denoted by $h\o m\m h\tr m, m\o h\m m\tl h, m\m m_{[-1]}\o m_{[0]}, m\m m_{(0)}\o m_{(1)}$, for any $h\in H$ and $m\in M$), such that $M$ is a left-left Yetter-Drinfel'd module, a left-right Long module, a right-right Yetter-Drinfel'd module and a right-left Long module, i.e.
\begin{align}
&(h_{1}\tr m)_{[-1]}h_{2}\o (h_{1}\tr m)_{[0]}=h_{1}m_{[-1]}\o h_{2}\tr m_{[0]},\label{y1}
\\&(h\tr m)_{(0)}\o (h\tr m)_{(1)}=h\tr m_{(0)}\o m_{(1)},\label{y2}
\\&(m\tl h_{2})_{(0)}\o h_{1}(m\tl h_{2})_{(1)}=m_{(0)}\tl h_{1} \o m_{(1)}h_{2},\label{y3}
\\&(m\tl h)_{[-1]}\o (m\tl h)_{[0]}=m_{[-1]}\o m_{[0]}\tl h.\label{y4}
\end{align}
\end{definition}
We denote by $\mathcal{LR}(H)$ the category whose objects are all Yetter-Drinfel'd-Long bimodule $M$ over $H$, the morphisms in the category are morphisms of $H$-bilinear and $H$-bicolinear.

If $H$ has a bijective antipode $S$, $\mathcal{LR}(H)$ becomes a strict braided monoidal category with the following structures: for any $M, N\in \mathcal{LR}(H)$, and $h\in H$, $m\in M$ and $n\in N$,
\begin{align*}
h\tr (m\o n)&=h_{1}\tr m\o h_{2}\tr n,
\\(m\o n)_{[-1]}\o (m\o n)_{[0]}&=m_{[-1]}n_{[-1]}\o m_{[0]}\o n_{[0]},
\\(m\o n)\tl h&=m\tl h_{1}\o n\tl h_{2},
\\(m\o n)_{(0)}\o (m\o n)_{(1)}&=m_{(0)}\o n_{(0)}\o m_{(1)}n_{(1)},
\end{align*}
the braiding
\begin{align*}
\psi_{M,N}:M\o N\ra N\o M : m\o n\m    m_{[-1]}\tr n_{(0)}\o m_{[0]}\tl n_{(1)}
\end{align*}
and the inverse
\begin{align*}
\psi_{N, M}^{-1}:N\o M\ra M\o N : n\o m\m m_{[0]}\tl S^{-1}(n_{(1)})\o S^{-1}(m_{[-1]})\tr n_{(0)}.
\end{align*}
\begin{definition}(\cite{M93})
A quasitriangular (QT) Hopf algebra is a pair $(H, R)$, where $H$ is a Hopf algebra over $\Bbbk$ and $R=R^{1}\o R^{2}\in H\o H$ is invertible, such that the following conditions hold ($r=R$):
\begin{enumerate}
\item[]
    \begin{enumerate}
\item[$(QT1)$]$\D(R^{1})\o R^{2}=R^{1}\o r^{1}\o R^{2}r^{2}$;
\item[$(QT2)$]$R^{1}\o \D(R^{2})=R^{1}r^{1}\o r^{2}\o R^{2}$;
\item[$(QT3)$] $\D^{cop}(h)R=R\D(h)$;
\item[$(QT4)$] $\v(R^{1})R^{2}=1=R^{1}\v(R^{2})$;
\item[$(QT5)$] If $R^{-1}=R^{2}\o R^{1}$, then $(H, R)$ is called a triangular Hopf algebra.
\end{enumerate}
\end{enumerate}
\end{definition}
\begin{definition}(\cite{M93})
A coquasitriangular (CQT) Hopf algebra is a pair $(H, \z)$, where $H$ is a Hopf algebra over $\Bbbk$ and $\z: H\o H\ra \Bbbk$ is a $\Bbbk$-bilinear form (braiding) which is convolution invertible in Hom $\!_{\Bbbk}(H\o H, \Bbbk)$ such that the following conditions hold:
   \begin{enumerate}
\item[]
\begin{enumerate}
\item[$(CQT1)$] $\z(h, gl)=\z(h_{1}, g)\z(h_{2},l)$;
\item[$(CQT2)$] $\z(hg, l)=\z(h , l_{2})\z(g,l_{1})$;
\item[$(CQT3)$] $\z(h_{1},g_{1})g_{2}h_{2}=h_{1}g_{1}\z(h_{2}, g_{2})$;
\item[$(CQT4)$] $\z(h, 1)=\v(h)=\z(1, h)$;
\item[$(CQT5)$] If $\z(h_{1}, g_{1})\z(g_{2}, h_{2})=\v(g)\v(h)$, then $(H, \z)$ is called a cotriangular Hopf algebra.
\end{enumerate}
\end{enumerate}
\end{definition}
The following are some examples of objects in $\mathcal{LR}(H)$.
\begin{example}\label{E2.1}
Let $H$ be a Hopf algebra. Then

$(1)$ $H_{1}=H\o H$ is a Yetter-Drinfel'd-Long bimodule with the following structures, for any $h, k, l\in H $:
\begin{align*}
&h\tr (k\o l)=hk\o l,&& \r^{l}(k\o l)=(k\o l)_{[-1]}\o (k\o l)_{[0]}=k_{1}S(k_{3})\o (k_{2}\o l),
\\&(k\o l)\tl h=k\o S(h_{1})l h_{2},&& \r^{r}(k\o l)=(k\o l)_{(0)}\o (k\o l)_{(1)}=(k\o l_{1})\o l_{2}.
\end{align*}

$(2)$ $H_{2}=H\o H$ is a Yetter-Drinfel'd-Long bimodule with the following structures, for any $h, k, l\in H $:
\begin{align*}
&h\tr (k\o l) =h_{1}kS(h_{2})\o l,&& \r^{l}(k\o l)=(k\o l)_{[-1]}\o (k\o l)_{[0]}=k_{1} \o (k_{2}\o l),
\\&(k\o l)\tl h= k\o lh,&& \r^{r}(k\o l)=(k\o l)_{(0)}\o (k\o l)_{(1)}=(k\o l_{2})\o S(l_{1})l_{3}.
\end{align*}

$(3)$ $H_{3}=H\o H$ is a Yetter-Drinfel'd-Long bimodule with the following structures, for any $h, k, l\in H $:
\begin{align*}
&h\tr (k\o l)= hk\o l,&& \r^{l}(k\o l)=(k\o l)_{[-1]}\o (k\o l)_{[0]}=k_{1}S(k_{3})\o (k_{2}\o l),
\\&(k\o l)\tl h= k\o lh,&& \r^{r}(k\o l)=(k\o l)_{(0)}\o (k\o l)_{(1)}=(k\o l_{2})\o S(l_{1})l_{3}.
\end{align*}

$(4)$ $H_{4}=H\o H$ is a Yetter-Drinfel'd-Long bimodule with the following structures, for any $h, k, l\in H $:
\begin{align*}
&h\tr (k\o l) =h_{1}kS(h_{2})\o l,&& \r^{l}(k\o l)=(k\o l)_{[-1]}\o (k\o l)_{[0]}=k_{1} \o (k_{2}\o l),
\\&(k\o l)\tl h=k\o S(h_{1})l h_{2},&& \r^{r}(k\o l)=(k\o l)_{(0)}\o (k\o l)_{(1)}=(k\o l_{1})\o l_{2}.
\end{align*}
Note that $H\o H$ is also a Hopf algebra with usual tensor product and usual tensor coproduct.
\end{example}
\section{Symmetric Yetter-Drinfel'd-Long categories}
\def\theequation{3.\arabic{equation}}
\setcounter{equation} {0}
In this section, we give necessary and sufficient conditions for Yetter-Drinfel'd-Long category $\mathcal{LR}(H)$ to be symmetric and pseudosymmetric, respectively.

Let $\mathcal{C}$ be a monoidal category and $\psi$ a braiding on $\mathcal{C}$. The braiding $\psi$ is called a symmetry if $\psi_{W, V}\circ \psi_{V, W}=id_{V\o W}$ for any $V, W\in \mathcal{C}$. In this case, $\mathcal{C}$ is called a symmetric braided category (see \cite{JS93}). The braiding $\psi$ is called a pseudosymmetry if the following condition holds, for any $U,V,W\in \mathcal{C}$:
\begin{align*}
(id_{W}\o \psi_{U,V})(\psi^{-1}_{W,U}\o id_{V})(id_{U}\o\psi_{V, W})=(\psi_{V,W}\o id_{U})(id_{V}\o\psi^{-1}_{W, U})(\psi_{U,V}\o id_{W}).
\end{align*}
In this case, $\mathcal{C}$ is called a pseudosymmetric braided category (see \cite{PSV10}).

Note that if $\psi$ is a symmetry, that is, $\psi^{-1}_{W, V}=\psi_{V, W}$, then obviously $\psi$ is a pseudosymmetry.
\begin{theorem}\label{T3.1}
Let $H$ be a Hopf algebra such that the canonical braiding of the Yetter-Drinfel'd-Long category $\mathcal{LR}(H)$ is a symmetry if and only if $H=\Bbbk$.
\end{theorem}
\begin{proof}
By Example \ref{E2.1}, $H_{1}$ and $H_{2}$ are two Yetter-Drinfel'd-Long bimodules. If the canonical braiding $\psi$ is a symmetry, that is, $\psi_{H_{2},H_{1}}\circ \psi_{H_{1},H_{2}}=id_{H_{1}\o H_{2}}$. Apply $\psi_{H_{2},H_{1}}\circ \psi_{H_{1},H_{2}}$ to the element $1\o k\o 1\o 1\in H_{1}\o H_{2}$, we have
\begin{align*}
&\psi_{H_{2},H_{1}}\circ \psi_{H_{1},H_{2}}(1\o k\o 1\o 1)
\\&=\psi_{H_{2},H_{1}}((1\o k)_{[-1]}\tr (1\o 1)_{(0)}\o (1\o k)_{[0]}\tl (1\o 1)_{(1)})
\\&=\psi_{H_{2},H_{1}}(1\tr (1\o 1) \o (1\o k) \tl 1)
\\&=\psi_{H_{2},H_{1}}(1\o 1 \o 1\o k)
\\&=(1\o 1)_{[-1]}\tr (1\o k)_{(0)}\o (1\o 1)_{[0]}\tl (1\o k)_{(1)}
\\&=1\tr (1\o k_{1}) \o (1\o 1) \tl k_{2}
\\&=1\o k_{1} \o 1\o k_{2}.
\end{align*}
Thus we have $1\o k\o 1\o 1=1\o k_{1} \o 1\o k_{2}$. Apply $\v\o \v\o \v \o id$ to both sides of the equation, we have $\v(k)1_{H}=k$. So $H=\Bbbk$.

The converse is straightforward, This completes the proof.
\end{proof}
Here, we will give three other proofs of Theorem \ref{T3.1}, and they are different from each other.
\begin{enumerate}
\item[$\bu$] By Example \ref{E2.1}, $H_{1}$ and $H_{3}$ are two Yetter-Drinfel'd-Long bimodules. If canonical braiding is a symmetry, that is, $\psi_{H_{3},H_{1}}\circ \psi_{H_{1},H_{3}}=id_{H_{1}\o H_{3}}$. For any $1\o k\o 1\o 1\in H_{1}\o H_{3}$, we easily get that $\psi_{H_{3},H_{1}}\circ \psi_{H_{1},H_{3}}(1\o k\o 1\o 1)=1\o k_{1} \o 1\o k_{2}$.

    Thus we have $1\o k\o 1\o 1=1\o k_{1} \o 1\o k_{2}$. Apply $\v\o \v\o \v \o id$ to both sides of the equation, we have $\v(k)1_{H}=k$. So $H=\Bbbk$.
\item[$\bu$] By Example \ref{E2.1}, $H_{2}$ and $H_{4}$ are two Yetter-Drinfel'd-Long bimodules. If canonical braiding is a symmetry, that is, $\psi_{H_{2},H_{4}}\circ \psi_{H_{4},H_{2}}=id_{H_{4}\o H_{2}}$. For any $1\o k\o 1\o 1\in H_{4}\o H_{2}$, we easily get that $\psi_{H_{2},H_{4}}\circ \psi_{H_{4},H_{2}}(1\o k\o 1\o 1)=1\o k_{1} \o 1\o k_{2}$.

    Thus we have $1\o k\o 1\o 1=1\o k_{1} \o 1\o k_{2}$. Apply $\v\o \v\o \v \o id$ to both sides of the equation, we have $\v(k)1_{H}=k$. So $H=\Bbbk$.
    \item[$\bu$] By Example \ref{E2.1}, $H_{3}$ and $H_{4}$ are two Yetter-Drinfel'd-Long bimodules. If canonical braiding is a symmetry, that is, $\psi_{H_{3},H_{4}}\circ \psi_{H_{4},H_{3}}=id_{H_{4}\o H_{3}}$. For any $1\o k\o 1\o 1\in H_{4}\o H_{3}$, we easily get that $\psi_{H_{3},H_{4}}\circ \psi_{H_{4},H_{3}}(1\o k\o 1\o 1)=1\o k_{1} \o 1\o k_{2}$.

    Thus we have $1\o k\o 1\o 1=1\o k_{1} \o 1\o k_{2}$. Apply $\v\o \v\o \v \o id$ to both sides of the equation, we have $\v(k)1_{H}=k$. So $H=\Bbbk$.
\end{enumerate}
If $H_{1}=\Bbbk\o H$ and $H_{2}=\Bbbk\o H$, then $H_{1}$ and $H_{2}$ are two right-right Yetter-Drinfel'd modules. Hence using Theorem \ref{T3.1}, we can improve the main result in \cite{P01}.
\begin{corollary}
Let $H$ be a Hopf algebra such that the canonical braiding of right-right Yetter-Drinfel'd category $\mathcal{YD}^{H}_{H}$ is a symmetry. Then $H=\Bbbk$.
\end{corollary}

In the following, we will introduce the pseudosymmetry on $\mathcal{LR}(H)$ over a Hopf algebra $H$. For this purpose, we need the following Lemma.
\begin{lemma}\label{L2.1}
Let $H$ be a cocommutative Hopf algebra. Then the canonical braiding $\psi_{H_{1},H_{2}}$ of the category $\mathcal{LR}(H)$ is the usual flip map.
\end{lemma}
\begin{proof}
For any $g\o h\o k\o l\in H_{1}\o H_{2}$, we have
\begin{align*}
\psi_{H_{1},H_{2}}&(g\o h\o k\o l)
\\&=(g\o h)_{[-1]}\tr (k\o l)_{(0)}\o (g\o h)_{[0]}\tl (k\o l)_{(1)}
\\&=g_{1}S(g_{3})\tr (k\o l_{2})\o (g_{2}\o h)\tl l_{1}S(l_{3})
\\&=g_{1}S(g_{2})\tr (k\o l_{3})\o (g_{3}\o h)\tl l_{1}S(l_{2}) \quad by ~cocommutative
\\&=1\tr (k\o l)\o (g\o h)\tl 1
\\&=k\o l\o g\o h.
\end{align*}
This completes the proof.
\end{proof}
We now give necessary and sufficient conditions for the canonical braiding of the category $\mathcal{LR}(H)$ to be a pseudosymmetry, we prove the necessary condition by a new method which is different from Proposition 2.5 in \cite{PM13}.
\begin{theorem}\label{T3.2}
Let $H$ be a Hopf algebra. Then the canonical braiding of the category $\mathcal{LR}(H)$ is pseudosymmetric if and only if $H$ is commutative and cocommutative.
\end{theorem}
\begin{proof}
Assume that the canonical braiding $\psi$ of the category $\mathcal{LR}(H)$ is pseudosymmetric.
We first check that $H$ is cocommutative.
For any $1\o 1\o k\o 1\o 1\o 1\in H_{1}\o H_{2}\o H_{1}$, we have
\begin{align*}
(id&\o \psi_{H_{1},H_{2}})\circ (\psi_{H_{1},H_{1}}^{-1}\o id)\circ(id\o \psi_{H_{2},H_{1}})(1\o 1\o k\o 1\o 1\o 1)
\\&=(id\o \psi_{H_{1},H_{2}})\circ (\psi_{H_{1},H_{1}}^{-1}\o id)(1\o 1\o (k\o 1)_{[-1]}\tr (1\o 1)_{(0)}\o (k\o 1)_{[0]}\tl (1\o 1)_{(1)})
\\&=(id\o \psi_{H_{1},H_{2}})\circ (\psi_{H_{1},H_{1}}^{-1}\o id)(1\o 1\o k_{1}\tr (1\o 1) \o (k_{2}\o 1)\tl 1)
\\&=(id\o \psi_{H_{1},H_{2}})\circ (\psi_{H_{1},H_{1}}^{-1}\o id)(1\o 1\o k_{1}\o 1 \o k_{2}\o 1)
\\&=(id\o \psi_{H_{1},H_{2}})((k_{1}\o 1)_{[0]}\tl S^{-1}((1\o 1)_{(1)})\o S^{-1}((k_{1}\o 1)_{[-1]})\tr (1\o 1)_{(0)} \o k_{2}\o 1)
\\&=(id\o \psi_{H_{1},H_{2}})((k_{2}\o 1)\tl 1\o S^{-1}(k_{1}S(k_{3}))\tr (1\o 1)  \o k_{4}\o 1)
\\&=(id\o \psi_{H_{1},H_{2}})(k_{2}\o 1\o k_{3}S^{-1}(k_{1})\o 1  \o k_{4}\o 1)
\\&=k_{2}\o 1\o (k_{3}S^{-1}(k_{1})\o 1)_{[-1]}\tr (k_{4}\o 1)_{(0)}\o (k_{3}S^{-1}(k_{1})\o 1)_{[0]}\tl (k_{4}\o 1)_{(1)}
\\&=k_{2}\o 1\o (k_{3}S^{-1}(k_{1}))_{1}S((k_{3}S^{-1}(k_{1}))_{3})\tr (k_{4}\o 1)\o ((k_{3}S^{-1}(k_{1}))_{2} \o 1)\tl 1
\\&=k_{2}\o 1\o [(k_{3}S^{-1}(k_{1}))_{1}S((k_{3}S^{-1}(k_{1}))_{3})]_{1}k_{4}S([(k_{3}S^{-1}(k_{1}))_{1}S((k_{3}S^{-1}(k_{1}))_{3})]_{2})\o 1
\\&\quad \o (k_{3}S^{-1}(k_{1}))_{2} \o 1
\end{align*}
and
\begin{align*}
(\psi_{H_{2},H_{1}}&\o id)\circ(id\o \psi_{H_{1},H_{1}}^{-1})\circ (\psi_{H_{1},H_{2}}\o id)(1\o 1\o k\o 1\o 1\o 1)
\\&=(\psi_{H_{2},H_{1}}\o id)\circ(id\o \psi_{H_{1},H_{1}}^{-1})
\\&\quad ((1\o 1)_{[-1]}\tr (k\o 1)_{(0)}\o (1\o 1)_{[0]}\tl (k\o 1)_{(1)}\o 1\o 1)
\\&=(\psi_{H_{2},H_{1}}\o id)\circ(id\o \psi_{H_{1},H_{1}}^{-1})(1\tr (k\o 1)\o (1\o 1)\tl 1\o 1\o 1)
\\&=(\psi_{H_{2},H_{1}}\o id)\circ(id\o \psi_{H_{1},H_{1}}^{-1})(k\o 1\o 1\o 1\o 1\o 1)
\\&=(\psi_{H_{2},H_{1}}\o id)(k\o 1\o 1\o 1\o 1\o 1)
\\&= (k\o 1)_{[-1]}\tr (1\o 1)_{(0)}\o (k\o 1)_{[0]}\tl (1\o 1)_{(1)}\o 1\o 1
\\&=k_{1}\tr (1\o 1)\o (k_{2}\o 1)\tl 1\o 1\o 1
\\&=k_{1}\o 1\o k_{2}\o 1\o 1\o 1.
\end{align*}
By assumption, $\mathcal{LR}(H)$ is pseudosymmetric, it follows that
\begin{align*}
k_{1}\o 1\o k_{2}\o 1\o 1\o 1= k_{2}&\o 1\o [(k_{3}S^{-1}(k_{1}))_{1}S((k_{3}S^{-1}(k_{1}))_{3})]_{1}k_{4}
\\&\times S([(k_{3}S^{-1}(k_{1}))_{1}S((k_{3}S^{-1}(k_{1}))_{3})]_{2})\o 1\o (k_{3}S^{-1}(k_{1}))_{2} \o 1
\end{align*}
Apply $id\o \v\o \v\o \v\o id\o \v$ to both sides of the above equation, we get $k_{2}\o k_{3}S^{-1}(k_{1})=k\o 1$. Therefore, we have
\begin{align*}
k_{2}\o k_{1}=k_{2}\o 1 k_{1}=k_{3}\o k_{4}S^{-1}(k_{2}) k_{1}=k_{1}\o k_{2}.
\end{align*}
So $H$ is cocommutative.

Next, we verify that $H$ is commutative. For any $1\o 1\o k\o 1\o g\o 1\in H_{1}\o H_{2}\o H_{2}$, we have
\begin{align*}
(id&\o \psi_{H_{1},H_{2}})\circ (\psi_{H_{2},H_{1}}^{-1}\o id)\circ(id\o \psi_{H_{2},H_{2}})(1\o 1\o k\o 1\o g\o 1)
\\&=(id\o \psi_{H_{1},H_{2}})\circ (\psi_{H_{2},H_{1}}^{-1}\o id)(1\o 1\o (k\o 1)_{[-1]}\tr (g\o 1)_{(0)}\o (k\o 1)_{[0]}\tl (g\o 1)_{(1)})
\\&=(id\o \psi_{H_{1},H_{2}})\circ (\psi_{H_{2},H_{1}}^{-1}\o id)(1\o 1\o k_{1}\tr (g\o 1) \o (k_{2}\o 1)\tl 1)
\\&=(id\o \psi_{H_{1},H_{2}})\circ (\psi_{H_{2},H_{1}}^{-1}\o id)(1\o 1\o k_{1}gS(k_{2})\o 1 \o k_{3}\o 1)
\\&=(id\o \psi_{H_{1},H_{2}})((k_{1}gS(k_{2})\o 1)_{[0]}\tl S^{-1}((1\o 1)_{(1)})
\\&\quad\o S^{-1}((k_{1}gS(k_{2})\o 1)_{[-1]})\tr (1\o 1)_{(0)} \o k_{3}\o 1)
\\&=(id\o \psi_{H_{1},H_{2}})((k_{2}g_{2}S(k_{3})\o 1)\tl 1\o S^{-1}(k_{1}g_{1}S(k_{4}))\tr (1\o 1)  \o k_{5}\o 1)
\\&=(id\o \psi_{H_{1},H_{2}})(k_{2}g_{2}S(k_{3})\o 1\o S^{-1}(k_{1}g_{1}S(k_{4}))\o 1  \o k_{5}\o 1)
\\&=k_{2}g_{2}S(k_{3})\o 1  \o k_{5}\o 1\o S^{-1}(k_{1}g_{1}S(k_{4}))\o 1 \quad by ~Lemma~\ref{L2.1}
\end{align*}
and
\begin{align*}
(\psi&_{H_{2},H_{2}}\o id)\circ(id\o \psi_{H_{2},H_{1}}^{-1})\circ (\psi_{H_{1},H_{2}}\o id)(1\o 1\o k\o 1\o g\o 1)
\\&=(\psi_{H_{2},H_{2}}\o id)\circ(id\o \psi_{H_{2},H_{1}}^{-1})(k\o 1 \o 1\o 1\o g\o 1) \quad by ~Lemma~\ref{L2.1}
\\&=(\psi_{H_{2},H_{2}}\o id)(k\o 1 \o (g\o 1)_{[0]}\tl S^{-1}((1\o 1)_{(1)})\o S^{-1}((g\o 1)_{[-1]})\tr (1\o 1)_{(0)})
\\&=(\psi_{H_{2},H_{2}}\o id)(k\o 1 \o (g_{2}\o 1) \tl 1 \o S^{-1}(g_{1})\tr (1\o 1) )
\\&=(\psi_{H_{2},H_{2}}\o id)(k\o 1 \o g_{2}\o 1 \o S^{-1}(g_{1})\o 1 )
\\&=(k\o 1)_{[-1]}\tr (g_{2}\o 1)_{(0)}\o (k\o 1)_{[0]}\tl (g_{2}\o 1)_{(1)} \o S^{-1}(g_{1})\o 1
\\&=k_{1}\tr (g_{2}\o 1) \o (k_{2}\o 1)\tl 1\o S^{-1}(g_{1})\o 1
\\&=k_{1}g_{2}S(k_{2})\o 1 \o k_{3}\o 1\o S^{-1}(g_{1})\o 1.
\end{align*}
Since $\mathcal{LR}(H)$ is pseudosymmetric, we get
\begin{align*}
k_{2}g_{2}S(k_{3})\o 1  \o k_{5}\o 1\o S^{-1}(k_{1}g_{1}S(k_{4}))\o 1=k_{1}g_{2}S(k_{2})\o 1 \o k_{3}\o 1\o S^{-1}(g_{1})\o 1.
\end{align*}
Apply $(\v\o \v\o id\o \v\o id\o \v)(id\o id\o id\o id\o S\o id)$ to both sides of the above equation, we get
$k_{3}\o k_{1}gS(k_{2})=k \o  g$. Hence, we have
\begin{align*}
gk=k_{1}gS(k_{2})k_{3}=k_{1}g\v(k_{2})=kg.
\end{align*}
So $H$ is commutative.

The proof of the converse can refer to Proposition 2.5 in \cite{PM13}. This completes the proof.
\end{proof}
If we consider $H_{1}=H\o \Bbbk$ and $H_{2}=H\o \Bbbk$, then $H_{1}$ and $H_{2}$ are two left-left Yetter-Drinfel'd modules. By the proof of Theorem \ref{T3.2}, we have the following result:
\begin{corollary}
The canonical braiding of $\!_{H}^{H}\mathcal{YD}$ is pseudosymmetric if and only if $H$ is cocommutative and commutative.
\end{corollary}
\section{The $u$-condition in $\mathcal{LR}(H)$}
\def\theequation{4.\arabic{equation}}
\setcounter{equation} {0}
In this section, we introduce the definition of the $u$-condition in $\mathcal{LR}(H)$ over Hopf algebra $H$ and discuss some properties and results related to the $u$-condition. It is easy to obtain the $u$-condition in $\!^{H}_{H}\mathcal{YD}$ when the right action and coaction are trivial.
\begin{definition}
Let $H$ be a Hopf algebra and $M\in \mathcal{LR}(H)$. Then $M$ is said to satisfy the $u$-condition if
\begin{align}\label{e3.1}
m_{[-1]}\tr m_{[0](0)}\tl m_{[0](1)}=m,
\end{align}
for any $m\in M$.
\end{definition}
Note that Eq.$(\ref{e3.1})$ is equivalent to the following equation:
\begin{align}\label{e3.2}
m_{(0)[-1]}\tr m_{(0)[0]}\tl m_{(1)}=m,
\end{align}
for any $m\in M$.

In the following, we will give a necessary and sufficient condition for $H_{1}$, $H_{2}$, $H_{3}$ and $H_{4}$ in Example \ref{E2.1} to satisfy the $u$-condition.
\begin{proposition}\label{P3.1}
Let $H$ be a Hopf algebra. Then
\begin{enumerate}
\item[$(1)$] $H_{1}$ satisfies the $u$-condition if and only if $S^{2}=id$.
\item[$(2)$] $H_{2}$ satisfies the $u$-condition if and only if $S^{2}=id$.
\item[$(3)$] $H_{3}$ satisfies the $u$-condition if and only if $S^{2}=id$.
\item[$(4)$] $H_{4}$ satisfies the $u$-condition if and only if $S^{2}=id$.
\end{enumerate}
\end{proposition}
\begin{proof}
It is basic in \cite{K95} that $S^{2}=id$ if and only if $S(h_{2})h_{1}=\v(h)$ or $h_{2}S(h_{1})=\v(h)$.

For $(1)$, if $S^{2}=id$, we only need to check that Eq.$(\ref{e3.1})$ holds. For any $k,l \in H$, we have
\begin{align*}
(k\o l)_{[-1]}&\tr (k\o l)_{[0](0)}\tl (k\o l)_{[0](1)}
\\&=k_{1}S(k_{3})\tr (k_{2}\o l)_{(0)}\tl (k_{2}\o l)_{(1)}
\\&=k_{1}S(k_{3})\tr (k_{2}\o l_{1})\tl l_{2}
\\&=k_{1}S(k_{3})k_{2}\o  S(l_{2})l_{1}l_{3}
\\&=k_{1}\v(k_{2})\o \v(l_{1})l_{2}
\\&=k\o l.
\end{align*}

Conversely, assume that $H_{1}$ satisfies the $u$-condition. For any $k\o 1\in H_{1}$, we have
\begin{align*}
(k\o 1)_{[-1]}&\tr (k\o 1)_{[0](0)}\tl (k\o 1)_{[0](1)}
\\&=k_{1}S(k_{3})\tr (k_{2}\o 1)_{(0)}\tl (k_{2}\o 1)_{(1)}
\\&=k_{1}S(k_{3})\tr (k_{2}\o 1)\tl 1
\\&=k_{1}S(k_{3})k_{2}\o  1.
\end{align*}
By assumption, we have $k_{1}S(k_{3})k_{2}\o  1=k\o 1$. Apply $id\o \v$ to both sides, we get
\begin{align}\label{e3.3}
k_{1}S(k_{3})k_{2}=k.
\end{align}
By computing we have
\begin{align*}
S(k_{2})k_{1}&=\v(k_{1})S(k_{3})k_{2}
\\&=(S(k_{1})k_{2})S(k_{4})k_{3}
\\&=S(k_{1})(k_{2}S(k_{4})k_{3})
\\&=S(k_{1})k_{2} \quad By~ Applying~ (\ref{e3.3})~ to~ k_{2}
\\&=\v(k).
\end{align*}
Hence $S^{2}=id$.

For $(2)$, if $S^{2}=id$, for any $k,l \in H$, we have
\begin{align*}
(k\o l)_{[-1]}&\tr (k\o l)_{[0](0)}\tl (k\o l)_{[0](1)}
\\&=k_{1}\tr (k_{2}\o l)_{(0)}\tl (k_{2}\o l)_{(1)}
\\&=k_{1}\tr (k_{2}\o l_{2})\tl  S(l_{1})l_{3}
\\&=k_{1}k_{3}S(k_{2})\o  l_{2}S(l_{1})l_{3}
\\&=k_{1}\v(k_{2})\o \v(l_{1})l_{2}
\\&=k\o l.
\end{align*}

Conversely, assume that $H_{2}$ satisfies the $u$-condition. For any $k\o 1\in H_{2}$, we have
\begin{align*}
(k\o 1)_{[-1]}&\tr (k\o 1)_{[0](0)}\tl (k\o 1)_{[0](1)}
\\&=k_{1}\tr (k_{2}\o 1)_{(0)}\tl (k_{2}\o 1)_{(1)}
\\&=k_{1}\tr (k_{2}\o 1)\tl  1
\\&=k_{1}k_{3}S(k_{2})\o  1.
\end{align*}
By assumption, we have $k_{1}k_{3}S(k_{2})\o  1=k\o 1$. Apply $id\o \v$ to both sides, we get
\begin{align}\label{e3.4}
k_{1}k_{3}S(k_{2})=k.
\end{align}
By computing we have
\begin{align*}
k_{2}S(k_{1})&=\v(k_{1})k_{3}S(k_{2})
\\&=(S(k_{1})k_{2})k_{4}S(k_{3})
\\&=S(k_{1})(k_{2}k_{4}S(k_{3}))
\\&=S(k_{1})k_{2} \quad By~ Applying~ (\ref{e3.4})~ to~ k_{2}
\\&=\v(k).
\end{align*}
Hence $S^{2}=id$.

Similarly, we can check that the statements $(3)$ and $(4)$ hold.
\end{proof}
\begin{proposition}\label{P3.2}
Let $H$ be a Hopf algebra and $S^{2}=id$, and assume that $M$ and $N$ satisfy the $u$-condition. Then $M\o N$ satisfies the $u$-condition if and only if $\psi_{M,N}$ is a symmetry.
\end{proposition}
\begin{proof}
For any $m\in M$ and $n\in N$, we have
\begin{align*}
(m\o &n)_{[-1]}\tr (m\o n)_{[0](0)}\tl (m\o n)_{[0](1)}
\\&=(m_{[-1]} n_{[-1]})\tr (m_{[0]}\o n_{[0]})_{(0)}\tl (m_{[0]}\o n_{[0]})_{(1)}
\\&=(m_{[-1]} n_{[-1]})\tr (m_{[0](0)}\o n_{[0](0)})\tl (m_{[0](1)}n_{[0](1)})
\\&=m_{[-1]} \tr [n_{[-1]}\tr (m_{[0](0)}\o n_{[0](0)})\tl m_{[0](1)}]\tl n_{[0](1)}
\\&=m_{[-1]}\tr [n_{ [-1]1}\tr (m_{[0] (0) }\tl m_{[0] (1)1})\o (n_{[-1]2}\tr n_{[0](0) })\tl m_{[0] (1)2}]\tl n_{[0] (1)}
\\&=m_{[-1]}\tr [n_{(0) [-1]1}\tr (m_{[0] (0) }\tl m_{[0] (1)1})\o (n_{(0) [-1]2}\tr n_{(0) [0] })\tl m_{[0] (1)2}]\tl n_{(1)}\quad by ~(\ref{e2.c})
\\&=m_{[-1]}\tr [n_{(0) [-1]1}(n_{(0)[-1]4}S(n_{(0) [-1]3}))\tr (m_{[0] (0) }\tl m_{[0] (1)3})
\\&\quad \o (n_{(0) [-1]2}\tr n_{(0) [0] })\tl (S(m_{[0] (1)2})m_{[0] (1)1})m_{[0] (1)4}]\tl n_{(1)}\quad by ~S^{2}=id
\\&=m_{[-1]}\tr [(n_{(0) [-1]11}n_{(0)[-1]2})S(n_{(0) [-1]13})\tr (m_{[0] (0) }\tl m_{[0] (1)22})
\\&\quad \o (n_{(0) [-1]12}\tr n_{(0) [0] })\tl S(m_{[0] (1)21})(m_{[0] (1)1}m_{[0] (1)23})]\tl n_{(1)}
\\&=m_{[-1]}\tr [(n_{(0) [-1]1}n_{(0) [0][-1]})S(n_{(0) [-1]3})\tr (m_{[0] (0)(0)}\tl m_{[0] (1)2})
\\&\quad \o (n_{(0) [-1]2}\tr n_{(0) [0][0]})\tl S(m_{[0] (1)1})(m_{[0] (0)(1)}m_{[0] (1)3})]\tl n_{(1)}
\\&=m_{[-1]}\tr [(n_{(0) [-1]1}\tr n_{(0) [0]})_{[-1]}n_{(0) [-1]2}S(n_{(0) [-1]3})\tr (m_{[0] (0)}\tl m_{[0] (1)3})_{(0)}
\\&\quad \o (n_{(0) [-1]1}\tr n_{(0) [0]})_{[0]}\tl S(m_{[0] (1)1})m_{[0] (1)2}(m_{[0] (0)}\tl m_{[0] (1)3})_{(1)}]\tl n_{(1)}~ by ~(\ref{y1}),~(\ref{y3})
\\&=m_{[-1]}\tr [(n_{(0) [-1]}\tr n_{(0) [0]})_{[-1]}\tr (m_{[0] (0)}\tl m_{[0] (1)})_{(0)}
\\&\quad\o (n_{(0) [-1]}\tr n_{(0) [0]})_{[0]}\tl (m_{[0] (0)}\tl m_{[0] (1)})_{(1)}]\tl n_{(1)}
\\&=m_{[-1]}\tr [\psi_{N, M}(n_{(0) [-1]}\tr n_{(0) [0]}\o m_{[0] (0)}\tl m_{[0] (1)})]\tl n_{(1)}
\\&=\psi_{N, M}(m_{[-1]}\tr [n_{(0) [-1]}\tr n_{(0) [0]}\o m_{[0] (0)}\tl m_{[0] (1)}]\tl n_{(1)})
\\&=\psi_{N, M}(m_{[-1]1}n_{(0) [-1]}\tr n_{(0) [0]}\tl n_{ (1)1}\o m_{ [-1]2}\tr m_{[0] (0)}\tl m_{[0] (1)} n_{(1)2})
\\&=\psi_{N, M}(m_{[-1]}n_{(0)(0)[-1]}\tr n_{(0)(0)[0]}\tl n_{(0)(1)}\o m_{[0][-1]}\tr m_{[0][0](0)}\tl m_{[0][0](1)} n_{(1)})
\\&=\psi_{N, M}(m_{[-1]}\tr n_{(0)}\o m_{[0]}\tl n_{(1)})\quad by ~(\ref{e3.1}),~(\ref{e3.2})
\\&=\psi_{N, M}\circ \psi_{M,N}(m\o n).
\end{align*}
This completes the proof.
\end{proof}
If we consider $M=H_{i}$ and $N=H_{j}$, for any $i, j=1,2,3,4$ (see Example $\ref{E2.1}$). By Proposition \ref{P3.1} and \ref{P3.2}, we obtain:
\begin{corollary}
Let $H$ be a Hopf algebra, and assume that $H_{i}$ and $H_{j}$ satisfy the $u$-condition. Then $H_{i}\o H_{j}$ satisfies the $u$-condition if and only if $\psi_{H_{i},H_{j}}$ is a symmetry, for any $i, j=1,2,3,4$.
\end{corollary}

\section{Yetter-Drinfel'd-Long categories over quasitriangular Hopf algebras}
\def\theequation{5.\arabic{equation}}
\setcounter{equation} {0}
In this section, we focus on $M\in \mathcal{LR}(H)$ for which $\psi_{M,M}$ is a symmetry. Triangular Hopf algebras give rise to such $M$.
\begin{theorem}\label{T5.1}
Let $(H, R)$ be a quasitriangular Hopf algebra. Then the category $\!_{H}\mathcal{M}_{H}$ of $H$-bimodules is a Yetter-Drinfel'd-Long subcategory of $\mathcal{LR}(H)$ under the coactions $\r^{l}(m)=R^{2}\o R^{1}\tr m$ and $\r^{r}(m)=m\tl R^{1}\o R^{2}$, where $\tr$ $(\tl, resp.)$ is the left (right, resp.) action on $M$.
\end{theorem}
\begin{proof}
First, we check that $M$ is a right $H$-comodule. By the definition of right $H$-comodule, for any $m\in M$, we have
\begin{align*}
(id\o \D)\r^{r}(m)&=(id\o \D)(m\tl R^{1}\o R^{2})
\\&=m\tl R^{1}\o R^{2}_{1}\o R^{2}_{2}
\\&=m\tl R^{1}r^{1}\o r^{2} \o R^{2}\quad by ~(QT2)
\\&=(\r^{r}\o id)(m\tl R^{1}\o R^{2})
\\&=(\r^{r}\o id)\r^{r}(m),
\end{align*}
and it is clear that $m_{(0)}\v(m_{(1)})=m\tl R^{1}\v(R^{2})=m\tl 1=m$. Similarly, we can get that $M$ is a left $H$-comodule.

Next, we verify the compatible condition of $H$-bicomodule. For any $m\in M$, we have
\begin{align*}
(id\o \r^{r})\r^{l}(m)&=(id\o \r^{r})(R^{2}\o R^{1}\tr m)
\\&=R^{2}\o (R^{1}\tr m)\tl r^{1}\o r^{2}
\\&=R^{2}\o R^{1}\tr (m\tl r^{1})\o r^{2}\quad by ~(\ref{e2.a})
\\&=(\r^{l}\o id)(m\tl r^{1}\o r^{2})
\\&=(\r^{l}\o id)\r^{r}(m).
\end{align*}

We now prove that $M$ satisfies the four compatible conditions $(\ref{y1})\sim (\ref{y4})$. Indeed, for any $h\in H$ and $m\in M$, we have
\begin{align*}
(h\tr m)_{(0)}\o (h\tr m)_{(1)}&=(h\tr m)\tl R^{1}\o R^{2}
\\&=h\tr (m\tl R^{1})\o R^{2}
\\&=h\tr m_{(0)}\o m_{(1)}.
\end{align*}
Thus Eq.$(\ref{y2})$ holds. For Eq.$(\ref{y3})$ , we have
\begin{align*}
m_{(0)}\tl h_{1}\o m_{(1)}h_{2}&=(m\tl  R^{1})\tl h_{1}\o R^{2}h_{2}
\\&=m\tl  R^{1}h_{1}\o R^{2}h_{2}
\\&=m\tl h_{2} R^{1}\o h_{1}R^{2}\quad by ~(QT3)
\\&=(m\tl h_{2})\tl R^{1}\o h_{1}R^{2}
\\&=(m\tl h_{2})_{(0)}\o h_{1}(m\tl h_{2})_{(1)}.
\end{align*}
Similarly, we can show that Eq.$(\ref{y1})$ and $(\ref{y4})$ hold.

Finally, we need to show that any morphisms in $\!_{H}\mathcal{M}_{H}$ are both left $H$-colinear and right $H$-colinear. For this purpose, we take any $M , N\in \!_{H}\mathcal{M}_{H}$, and assume that $f:M\ra N$ is a morphism in $\!_{H}\mathcal{M}_{H}$, we get
\begin{align*}
(f\o id)\circ\r^{r}_{M}(m)=f (m\tl R^{1})\o R^{2}=f (m)\tl R^{1}\o R^{2}=\r^{r}_{N}\circ f(m).
\end{align*}
So $f$ is right $H$-colinear. Similarly, we can obtain that $f$ described above is left $H$-colinear.

This completes the proof.
\end{proof}
\begin{proposition}\label{P5.1}
Let $H$ be a triangular Hopf algebra. Then the Yetter-Drinfel'd-Long subcategory $\!_{H}\mathcal{M}_{H}$ defined above is symmetric.
\end{proposition}
\begin{proof}
For any $m\in M$ and $n\in N$, we have
\begin{align*}
\psi_{N, M}\circ &\psi_{M, N}(m\o n)
\\&=\psi_{N, M}(R^{2}\tr n\tl r^{1}\o R^{1}\tr m\tl r^{2})
\\&=Q^{2}\tr (R^{1}\tr m\tl r^{2})\tl q^{1}\o Q^{1}\tr (R^{2}\tr n\tl r^{1})\tl q^{2}
\\&=Q^{2} R^{1}\tr m\tl r^{2} q^{1}\o Q^{1} R^{2}\tr n\tl r^{1} q^{2} \quad by ~ (QT5)
\\&=1\tr m\tl 1\o 1\tr n\tl 1
\\&=m\o n.
\end{align*}
Thus the subcategory $\!_{H}\mathcal{M}_{H}$ is symmetric.
\end{proof}
By Theorem \ref{T5.1} and Proposition \ref{P5.1}, we know that If $(H, R)$ be a triangular Hopf algebra then the subcategory $\!_{H}\mathcal{M}_{H}$ described above is symmetric. A particular example is $M=H\o H$.
In the following we prove the converse. That is, assume that the braiding $\psi_{H\o H, H\o H}$ is a symmetry forces $(H, R)$ to be triangular, where $H\o H$ is a Hopf algebra with usual tensor product and tensor coproduct.
\begin{theorem}\label{T4.3}
Let $H$ be a Hopf algebra with a bijective antipode, and assume that $(H\o H, \tr=m\o id, \r^{l}=\r_{1}\o id, \tl=id\o m, \r^{r}=id\o \r_{2})\in \mathcal{LR}(H)$, where $m$ is usual multiplication and $\r_{1}$ ($\r_{2}$, resp.) is a left (right, resp.) coaction on $H$. Then $\psi_{H\o H, H\o H}$ is a symmetry if and only if there exists $R\in H\o H$ so that $(H, R)$ is triangular. And then $\r^{l}$ and $\r^{r}$ are induced by $R$. That is,
\begin{align*}
&\r^{l}(k\o l)=R^{2}\o R^{1}k\o l,
\\&\r^{r}(k\o l)=k\o lR^{1}\o R^{2},
\end{align*}
for any $k, l\in H$, in particular, $R^{\tau}\o 1=\r^{l}(1\o 1)$ and $1\o R=\r^{r}(1\o 1)$.
\end{theorem}
\begin{proof}
If $\psi=\psi_{H\o H, H\o H}$ is a symmetry, for any $k, l, g, h\in H$, we have
\begin{align}
\nonumber\psi(k\o l\o g\o h)&= (k\o l)_{[-1]}\tr (g\o h)_{(0)}\o (k\o l)_{[0]}\tl (g\o h)_{(1)}
\\&=(g\o h)_{[0]}\tl S^{-1}((k\o l)_{(1)})\o S^{-1}((g\o h)_{[-1]})\tr (k\o l)_{(0)}. \label{e5.0}
\end{align}
In particular, let $\r^{l}(1\o 1)=x_{i}\o y_{i}\o 1$ and $\r^{r}(1\o 1)=1\o s_{i}\o t_{i}$. Then
\begin{align*}
x_{i}\o &s_{i}\o y_{i}\o t_{i}
\\&=x_{i}\tr (1\o s_{i})\o (y_{i}\o 1)\tl t_{i}
\\&=(1\o 1)_{[-1]}\tr (1\o 1)_{(0)}\o (1\o 1)_{[0]}\tl (1\o 1)_{(1)}
\\&=(1\o 1)_{[0]}\tl S^{-1}((1\o 1)_{(1)})\o S^{-1}((1\o 1)_{[-1]})\tr (1\o 1)_{(0)}\quad by ~(\ref{e5.0})
\\&=(y_{i}\o 1)\tl S^{-1}(t_{i})\o S^{-1}(x_{i})\tr (1\o s_{i})
\\&=y_{i}\o S^{-1}(t_{i})\o S^{-1}(x_{i})\o s_{i}.
\end{align*}
Thus
\begin{align*}
x_{i}\o s_{i}\o y_{i}\o t_{i}=y_{i}\o S^{-1}(t_{i})\o S^{-1}(x_{i})\o s_{i}.
\end{align*}
Apply $id\o \v\o id \o \v$ and $\v\o id \o \v\o id$ to both sides, respectively, we have
\begin{align}
&x_{i}\o y_{i}=y_{i}\o  S^{-1}(x_{i}),\label{e5.1}
\\&s_{i}\o t_{i}=S^{-1}(t_{i})\o s_{i}.\label{e5.2}
\end{align}
Apply $id\o S$ to Eq.$(\ref{e5.1})$ yields
\begin{align}
x_{i}\o S(y_{i})=y_{i}\o  x_{i}.\label{e5.3}
\end{align}

Set $R\o 1=y_{i}\o x_{i}\o  1=(\tau\o id)\circ\r^{l}(1\o 1)$ and $1\o R=1\o s_{i}\o t_{i}=\r^{r}(1\o 1)$.
In the following, we wish to show that $(H, R)$ is triangular and that $\r^{l}$ and $\r^{r}$ are induced by $R$. For this purpose, we first need the following equations
$\r^{l}(k\o l)=(id\o \v\o id^{2})\psi(k\o l\o 1\o 1)$ and $\r^{r}(k\o l)=(id^{2}\o \v\o id)\psi(1\o 1\o k\o l)$. Indeed, for any $k, l\in H$:
\begin{align*}
(id\o \v&\o id^{2})\psi(k\o l\o 1\o 1)
\\&=(id\o \v\o id^{2})((k\o l)_{[-1]}\tr (1\o 1)_{(0)}\o (k\o l)_{[0]}\tl (1\o 1)_{(1)})
\\&=(id\o \v\o id^{2})((k\o l)_{[-1]}\tr (1\o s_{i})\o (k\o l)_{[0]}\tl t_{i})
\\&=(id\o \v\o id^{2})((k\o l)_{[-1]}\o s_{i}\o (k\o l)_{[0]}\tl t_{i})
\\&=(id\o \v\o id^{2})((k\o l)_{[-1]}\o S^{-1}(t_{i})\o (k\o l)_{[0]}\tl s_{i})   \quad by ~(\ref{e5.2})
\\&=(k\o l)_{[-1]}\o  (k\o l)_{[0]}\tl 1
\\&=(k\o l)_{[-1]}\o  (k\o l)_{[0]}
\\&=\r^{l}(k\o l)
\end{align*}
and
\begin{align*}
(id^{2}&\o \v\o id)\psi(1\o 1\o k\o l)
\\&=(id^{2}\o \v\o id)((1\o 1)_{[-1]}\tr (k\o l)_{(0)}\o (1\o 1)_{[0]}\tl (k\o l)_{(1)})
\\&=(id^{2}\o \v\o id)(x_{i}\tr (k\o l)_{(0)}\o (y_{i}\o 1)\tl (k\o l)_{(1)})
\\&=(id^{2}\o \v\o id)(x_{i}\tr (k\o l)_{(0)}\o y_{i}\o (k\o l)_{(1)})
\\&=(id^{2}\o \v\o id)(y_{i}\tr (k\o l)_{(0)}\o S^{-1}(x_{i})\o (k\o l)_{(1)})  \quad by ~(\ref{e5.1})
\\&=1\tr (k\o l)_{(0)}\o  (k\o l)_{(1)}
\\&=(k\o l)_{(0)}\o  (k\o l)_{(1)}
\\&=\r^{r}(k\o l).
\end{align*}
We now prove that $\r^{l}$ and $\r^{r}$ are induced by $R$. For any $k, l\in H$, we have
\begin{align*}
\r^{l}&(k\o l)=(id\o \v\o id^{2})\psi(k\o l\o 1\o 1)
\\&=(id\o \v\o id^{2})((1\o 1)_{[0]}\tl S^{-1}((k\o l)_{(1)})\o S^{-1}((1\o 1)_{[-1]})\tr (k\o l)_{(0)}) \quad by ~(\ref{e5.0})
\\&=(id\o \v\o id^{2})((y_{i}\o 1)\tl S^{-1}((k\o l)_{(1)})\o S^{-1}(x_{i})\tr (k\o l)_{(0)})
\\&=(id\o \v\o id^{2})(y_{i}\o S^{-1}((k\o l)_{(1)})\o S^{-1}(x_{i})\tr (k\o l)_{(0)})
\\&=y_{i}\o  S^{-1}(x_{i})\tr (k\o l)
\\&=y_{i}\o  S^{-1}(x_{i})k\o l
\\&=x_{i}\o  y_{i}k\o l. \quad by ~(\ref{e5.1})
\end{align*}
and
\begin{align*}
\r^{r}&(k\o l)=(id^{2}\o \v\o id)\psi(1\o 1\o k\o l)
\\&=(id^{2}\o \v\o id)((k\o l)_{[0]}\tl S^{-1}((1\o 1)_{(1)})\o S^{-1}((k\o l)_{[-1]})\tr (1\o 1)_{(0)}) \quad by ~(\ref{e5.0})
\\&=(id^{2}\o \v\o id)((k\o l)_{[0]}\tl S^{-1}(t_{i})\o S^{-1}((k\o l)_{[-1]})\tr (1\o s_{i}))
\\&=(id^{2}\o \v\o id)((k\o l)_{[0]}\tl S^{-1}(t_{i})\o S^{-1}((k\o l)_{[-1]})\o s_{i})
\\&=(k\o l)\tl S^{-1}(t_{i})\o  s_{i}
\\&=k\o l S^{-1}(t_{i})\o  s_{i}
\\&=k\o l s_{i}\o  t_{i}. \quad by ~(\ref{e5.2})
\end{align*}
Thus
\begin{align}
&\r^{l}(k\o l)=x_{i}\o  y_{i}k\o l, \label{e5.4}
\\&\r^{r}(k\o l)=k\o l s_{i}\o  t_{i}. \label{e5.5}
\end{align}

Finally, we verify that $(H, R)$ is triangular. By definition, we need to prove the five equations (QT1) $\sim$ (QT5). For (QT1), we only have to check that $\D(y_{i})\o x_{i}=y_{i}\o y_{j}\o x_{i}x_{j}$.
\begin{align*}
\D(y_{i})&\o x_{i}=(id^{3} \o \v)(\D(y_{i})\o x_{i}\o 1)
\\&=(id^{3} \o \v)(\D(x_{i})\o S(y_{i})\o 1)  \quad by ~(\ref{e5.3})
\\&=(id^{2}\o S\o \v)(\D\o id^{2})(x_{i}\o y_{i}\o 1)
\\&=(id^{2}\o S\o \v)(\D\o id^{2})\r^{l}(1\o 1)
\\&=(id^{2}\o S\o \v)(id\o \r^{l})\r^{l}(1\o 1)
\\&=(id^{2}\o S\o \v)(x_{i}\o \r^{l}(y_{i}\o 1))
\\&=(id^{2}\o S\o \v)(x_{i}\o x_{j}\o y_{j}y_{i}\o 1))   \quad by ~(\ref{e5.4})
\\&=(id^{2}\o S\o \v)(y_{i}\o y_{j}\o S^{-1}(x_{j})S^{-1}(x_{i})\o 1))\quad by ~(\ref{e5.1})
\\&=y_{i}\o y_{j}\o x_{i}x_{j}.
\end{align*}
Similarly, we can check that (QT2) holds. For (QT3), we only need to show that $ h_{2}y_{i}\o h_{1}x_{i}=  y_{i}h_{1}\o x_{i}h_{2}$. Since both $\psi$ and $\v$ are $H$-module maps, we have
\begin{align*}
h_{1}x_{i}&\o h_{2}y_{i}=(id\o \v\o id\o \v)(h_{1}x_{i}\o1\o  h_{2}y_{i}\o 1)
\\&=(id\o \v\o id\o \v)(h_{1}\tr (x_{i}\o1)\o  h_{2}\tr (y_{i}\o 1))
\\&=(id\o \v\o id\o \v)[h\tr (x_{i}\o1\o y_{i}\o 1)]
\\&=h\tr [(id\o \v\o id\o \v)(x_{i}\o1\o y_{i}\o 1)]
\\&=h\tr [(id\o id\o \v)\circ \r^{l}(1\o 1)]
\\&=h\tr [(id\o \v\o id\o \v)\psi(1\o 1\o 1\o 1)]
\\&=(id\o \v\o id\o \v)[h\tr \psi(1\o 1\o 1\o 1)]
\\&=(id\o \v\o id\o \v)[\psi(h\tr (1\o 1\o 1\o 1))]
\\&=(id\o \v\o id\o \v)[\psi(h_{1}\o 1\o h_{2}\o 1)]
\\&=(id\o \v\o id\o \v)[(h_{1}\o 1)_{[-1]}\tr (h_{2}\o 1)_{(0)}\o (h_{1}\o 1)_{[0]}\tl (h_{2}\o 1)_{(1)}]
\\&=(id\o \v\o id\o \v)[x_{i}\tr (h_{2}\o s_{i})\o (y_{i}h_{1}\o 1)\tl t_{i}]\quad by ~(\ref{e5.4}),~(\ref{e5.5})
\\&=(id\o \v\o id\o \v)[x_{i}h_{2}\o s_{i}\o y_{i}h_{1}\o  t_{i}]
\\&=x_{i}h_{2}\o  y_{i}h_{1}.
\end{align*}
For (QT4), we have
\begin{align*}
\v(R^{1})R^{2}&=(\v\o id\o \v)(R^{1}\o R^{2}\o 1)
\\&=(\v\o id\o \v)(y_{i}\o x_{i}\o 1)
\\&=(\v\o id\o \v)(S^{-1}(x_{i})\o y_{i}\o 1) \quad by ~(\ref{e5.1})
\\&=(\v\o id\o \v)(x_{i}\o y_{i}\o 1)
\\&=(\v\o id\o \v)\r^{l}(1\o 1)
\\&=1.
\end{align*}
Similarly, we can check that $\v(R^{2})R^{1}=1$. For (QT5), we have
\begin{align*}
1\o 1&\o 1\o 1=\psi^{2}(1\o 1\o 1\o 1)
\\&=\psi((1\o 1)_{[-1]}\tr (1\o 1)_{(0)}\o (1\o 1)_{[0]}\tl (1\o 1)_{(1)})
\\&=\psi(x_{i}\tr (1\o s_{i})\o (y_{i}\o 1)\tl t_{i})
\\&=\psi(x_{i}\o s_{i}\o y_{i}\o t_{i})
\\&=(x_{i}\o s_{i})_{[-1]}\tr (y_{i}\o t_{i})_{(0)}\o (x_{i}\o s_{i})_{[0]}\tl (y_{i}\o t_{i})_{(1)}
\\&=x_{j}\tr (y_{i}\o t_{i}s_{j})\o (y_{j}x_{i}\o s_{i})\tl t_{j}
\\&=x_{j}y_{i}\o t_{i}s_{j}\o y_{j}x_{i}\o s_{i} t_{j}.
\end{align*}
Thus, $R$ is invertible and $R^{-1}=x_{i}\o y_{i}=t_{i}\o s_{i}$.

The converse is Theorem $\ref{T5.1}$ and Proposition $\ref{P5.1}$. This completes the proof.
\end{proof}
As a corollary we have:
\begin{corollary}
Let $H$ be a Hopf algebra with a bijective antipode. Then, for $H_{3}\in \mathcal{LR}(H)$, the braiding $\psi_{H_{3}, H_{3}}$ is a symmetry if and only if $H$ is cocommutative.
\end{corollary}
\begin{proof}
If the braiding satisfies $\psi_{H_{3}, H_{3}}^{2}=id$, then by Theorem \ref{T4.3} $(H, R)$ is triangular with $\r^{l}(1\o 1)=R^{\tau}\o 1$. Since $\r^{l}(k\o l)=k_{1}S(k_{3})\o k_{2}\o l$ for any $k,l\in H$, we have $\r^{l}(1\o 1)=1\o 1\o 1$, so $R=1\o 1$. Thus (QT3) implies that $H$ is cocommutative.

Conversely, assume that $H$ is cocommutative,
for any $k,l, g, h\in H$, we have
\begin{align*}
&\psi_{H_{3}, H_{3}}(k\o l\o g\o h)
\\&=(k\o l)_{[-1]}\tr (g\o h)_{(0)}\o (k\o l)_{[0]}\tl (g\o h)_{(1)}
\\&=k_{1}S(k_{3})\tr (g\o h_{2})\o (k_{2}\o l)\tl h_{1}S(h_{3})
\\&=k_{1}S(k_{2})\tr (g\o h_{3})\o (k_{3}\o l)\tl h_{1}S(h_{2})\quad by ~ H~is ~cocommutative
\\&=1\tr (g\o h) \o (k \o l) \tl 1
\\&=g\o h \o k\o l.
\end{align*}
It is clear that the braiding $\psi_{H_{3}, H_{3}}$ is a symmetry.
\end{proof}
If we consider $H\o \Bbbk$, by Theorem \ref{T4.3}, we generalize the important result in \cite{CW98}.
\begin{corollary}
Let $H$ be a Hopf algebra with a bijective antipode, and assume that $(H, m, \r )\in \!^{H}_{H}\mathcal{YD}$, where $m$ is usual multiplication. Then $\psi_{H,  H}$ is a symmetry if and only if there exists $R\in H\o H$ so that $(H, R)$ is triangular. And then $\r $ is induced by $R$. That is,
\begin{align*}
\r(k)=R^{2}\o R^{1}k,
\end{align*}
for any $k\in H$, in particular, $R^{\tau} =\r (1 )$.
\end{corollary}
\section{Yetter-Drinfel'd-Long categories over coquasitriangular Hopf algebras}
\def\theequation{6.\arabic{equation}}
\setcounter{equation} {0}
In this section, we discuss the dual cases of section 4.
\begin{theorem}\label{T6.1}
Let $(H, \z)$ be a coquasitriangular Hopf algebra. Then the category $\!^{H}\mathcal{M}^{H}$ of $H$-bicomodules is a Yetter-Drinfel'd-Long subcategory of $\mathcal{LR}(H)$ under the actions $h\tr m=\z(h, m_{[-1]})m_{[0]}$ and $m\tl h=m_{(0)}\z(h, m_{(1)})$, for any $h\in H$ and $m\in M\in \!^{H}\mathcal{M}^{H}$.
\end{theorem}
\begin{proof}
First, we prove that $(M, \tl)$ is a right $H$-module. For any $h, g\in H$ and $m\in M$, we have
\begin{align*}
(m\tl g )\tl h&=m_{(0)}\tl h \z(g, m_{(1)})
\\&=m_{(0)(0)}\z(h, m_{(0)(1)})\z(g, m_{(1)})
\\&=m_{(0)}\z(h, m_{(1)1})\z(g, m_{(1)2})
\\&=m_{(0)}\z(gh, m_{(1)})\quad by ~ (CQT2)
\\&=m\tl gh,
\end{align*}
and it is clear that $m\tl 1=m_{(0)}\z(1, m_{(1)})=m_{(0)}\v(m_{(1)})=m$. Similarly, we can obtain that $(M, \tr)$ is a left $H$-module.

Next, we check the compatible condition of $H$-bimodule. For any $h, g\in H$ and $m\in M$, we have
\begin{align*}
(h\tr m)\tl g&=\z(h, m_{[-1]})m_{[0]}\tl g
\\&=\z(h, m_{[-1]})m_{[0](0)}\z(g, m_{[0](1)})
\\&=\z(h, m_{(0)[-1]})m_{(0)[0]}\z(g, m_{(1)})\quad by ~(\ref{e2.c})
\\&=h\tr m_{(0)}\z(g, m_{(1)})
\\&=h\tr (m\tl g).
\end{align*}

We now check that the four compatible conditions $(\ref{y1})\sim (\ref{y4})$. For any $h\in H$ and $m\in M$, we have
\begin{align*}
(h\tr m&)_{(0)}\o (h\tr m)_{(1)}
\\&=\z(h, m_{[-1]})m_{[0](0)}\o (h\tr m)_{[0](1)}
\\&=\z(h, m_{(0)[-1]})m_{(0)[0]}\o m_{(1)} \quad by ~(\ref{e2.c})
\\&=h\tr m_{(0)}\o m_{(1)}.
\end{align*}
Thus Eq.$(\ref{y2})$ holds. For Eq.$(\ref{y3})$, we have
\begin{align*}
m_{(0)}&\tl h_{1}\o m_{(1)}h_{2}
\\&=m_{(0)(0)}\z(h_{1}, m_{(0)(1)})\o m_{(1)}h_{2}
\\&=m_{(0)}\o \z(h_{1}, m_{(1)1})m_{(1)2}h_{2}
\\&=m_{(0)}\o h_{1}m_{(1)1}\z(h_{2}, m_{(1)2})\quad by ~(CQT3)
\\&=m_{(0)(0)}\z(h_{2}, m_{(1)})\o h_{1}m_{(0)(1)}
\\&=(m\tl h_{2})_{(0)}\o h_{1}(m\tl h_{2})_{(1)}.
\end{align*}
Similarly, we can verify that Eq.$(\ref{y1})$ and $(\ref{y4})$ hold.

Finally, we have to prove that any morphisms in $\!^{H}\mathcal{M}^{H}$ are both left $H$-linear and right $H$-linear. For this purpose, we take any $M, N\in \!^{H}\mathcal{M}^{H}$, and assume that $f:M\ra N$ is a morphism in $\!^{H}\mathcal{M}^{H}$, we have
\begin{align*}
f(m \tl h)=f(m_{(0)})\z(h, m_{(1)})=f(m)_{(0)}\z(h, f(m)_{(1)})=f(m) \tl h.
\end{align*}
So $f$ is right $H$-linear. Similarly, we can obtain that $f$ is left $H$-linear.

This completes the proof.
\end{proof}
\begin{proposition}\label{P6.2}
Let $H$ be a cotriangular Hopf algebra. Then the Yetter-Drinfel'd-Long subcategory $\!^{H}\mathcal{M}^{H}$ defined above is symmetric.
\end{proposition}
\begin{proof}
For any $m\in M$ and $n\in N$, we have
\begin{align*}
\psi_{N,M}&\circ \psi_{M,N}(m\o n)
\\&=\psi_{N,M}(m_{[-1]}\tr n_{(0)}\o m_{[0]}\tl n_{(1)})
\\&=\psi_{N,M}(\z(m_{[-1]}, n_{(0)[-1]})n_{(0)[0]}\o m_{[0](0)}\z(n_{(1)}, m_{[0](1)}))
\\&=\z(m_{[-1]}, n_{(0)[-1]})\z(n_{(1)}, m_{[0](1)})n_{(0)[0][-1]}\tr m_{[0](0)(0)}\o n_{(0)[0][0]}\tl m_{[0](0)(1)}
\\&=\z(m_{[-1]}, n_{(0)[-1]1})\z(n_{(1)}, m_{[0](1)2})n_{(0)[-1]2}\tr m_{[0](0)}\o n_{(0) [0]}\tl m_{[0](1)1}
\\&=\z(m_{(0)[-1]}, n_{[-1]1})\z(n_{[0](1)}, m_{(1)2})n_{[-1]2}\tr m_{(0)[0]}\o n_{[0](0) }\tl m_{(1)1}\quad by ~(\ref{e2.c})
\\&=\z(m_{(0)[-1]}, n_{[-1]1})\z(n_{[0](1)}, m_{(1)2})
\\&\qquad\z(n_{[-1]2}, m_{(0)[0][-1]})m_{(0)[0][0]}\o n_{[0](0)(0)}\z(m_{(1)1}, n_{[0](0)(1)})
\\&=\z(m_{(0)[-1]1}, n_{[-1]1})\z(n_{[-1]2}, m_{(0)[-1]2})
\\&\qquad\z(m_{(1)1}, n_{[0](1)1})\z(n_{[0](1)2}, m_{(1)2})m_{(0)[0] }\o n_{[0](0)}\quad by ~(CQT5)
\\&=m\o n.
\end{align*}
So the subcategory $\!^{H}\mathcal{M}^{H}$ is symmetric.
\end{proof}
\begin{theorem}\label{T6.3}
Let $H$ be a Hopf algebra with a bijective antipode, and assume that $(H\o H, \tr=\rh\o id, \r^{l}=\D\o id, \tl=id\o \lh, \r^{r}=id\o \D)\in \mathcal{LR}(H)$, where $\D$ is usual comultiplication and $\rh$ ($\lh$, resp.) is a left (right, resp.) action on $H$. Then $\psi_{H\o H, H\o H}$ is a symmetry if and only if there exists a braiding $\z:H\o H\ra \Bbbk$ so that $(H, \z)$ is cotriangular Hopf algebra. And then $\z(k, g)\z(h, l)=(\v\o\v\o\v\o\v )\psi(k\o l\o g\o h)$, for any $k,l,g,h \in H$. That is,
\begin{align*}
&h\tr (k\o l)=h\rh k\o l=\z(h, k_{1})k_{2}\o l,
\\&(k\o l)\tl h=k\o l\lh h=k\o l_{1}\z(h, l_{2}).
\end{align*}
\end{theorem}
\begin{proof}
Assume that $\psi=\psi_{H\o H,H\o H}$ is a symmetry, then for any $k,l,g,h \in H$,
\begin{align*}
\psi(k\o l\o g\o h)&= (k\o l)_{[-1]}\tr (g\o h)_{(0)}\o (k\o l)_{[0]}\tl (g\o h)_{(1)}
\\&=(g\o h)_{[0]}\tl S^{-1}((k\o l)_{(1)})\o S^{-1}((g\o h)_{[-1]})\tr (k\o l)_{(0)},
\end{align*}
i.e.
\begin{align}
\nonumber \psi(k\o l\o g\o h)&=k_{1}\rh g\o h_{1}\o k_{2}\o l\lh h_{2}
\\&\quad=g_{2}\o h\lh S^{-1}(l_{2})\o S^{-1}(g_{1})\rh k\o l_{1}. \label{e6.3}
\end{align}
Define for any $k,l,g,h \in H$, $\z(k, g)\z(h, l)=(\v\o\v\o\v\o\v )\psi(k\o l\o g\o h)$.
Let $l=h=1$, and apply $\v\o \v\o \v\o \v$ to Eq.$(\ref{e6.3})$, we get
\begin{align}\label{e6.4}
\z(k, g)=\v(k\rh g)=\v(S^{-1}(g)\rh k)=\z(S^{-1}(g), k).
\end{align}
By applying $\z(k, g)=\z(S^{-1}(g), k)$ to $\z(g, S(k))$, we get
\begin{align}\label{e6.5}
\z(k, g)=\z(g, S(k)).
\end{align}
Similarly, we can get that
\begin{align}\label{e6.6}
\z(h, l)=\v(l\lh h)=\v(h\lh S^{-1}(l))=\z(S^{-1}(l), h)=\z(l, S(h)).
\end{align}
Moreover, let $l=h=1$, and apply $id\o \v\o \v\o \v$ to Eq.$(\ref{e6.3})$, we get by $(\ref{e6.4})$, that for any $k, g\in H$,
\begin{align}\label{e6.7}
k\rh g=\z(S^{-1}(g_{1}), k) g_{2}=\z(k, g_{1}) g_{2}.
\end{align}
Similarly, we can get by $(\ref{e6.6})$, that for any $l, h\in H$,
\begin{align*}
l\lh h=\z(S^{-1}(l_{2}), h) l_{1}=\z(h, l_{2})l_{1}.
\end{align*}
Thus we have
\begin{align*}
&h\tr (k\o l)=h\rh k\o l=\z(h, k_{1})k_{2}\o l,
\\&(k\o l)\tl h=k\o l\lh h=k\o l_{1}\z(h, l_{2}).
\end{align*}
By definition of cotriangular, we need to prove the five equations (CQT1) $\sim$ (CQT5). We prove (CQT2) first. For any $h, g,l\in H$, we have
\begin{align*}
\z(hg, l)&=\v(hg \rh l)
\\&=(\v\o \v\o \v\o\v)(h_{1}\rh (g \rh l)\o 1\o h_{2}\o 1)
\\&=(\v\o \v\o \v\o\v)(h_{1}\tr (g \rh l\o 1)\o (h_{2}\o 1)\tl 1)
\\&=(\v\o \v\o \v\o\v)\psi(h\o 1\o g \rh l\o 1)
\\&=\z(h, g \rh l)\z(1, 1)
\\&=\z(h, \z(g, l_{1}) l_{2})\quad by ~(\ref{e6.7})
\\&=\z(h,  l_{2})\z(g, l_{1}).
\end{align*}
Next we prove (CQT1). For any $h, g,l\in H$, we have
\begin{align*}
\z(h, gl)&=\z(gl, S(h)) \quad by ~(\ref{e6.5})
\\&=\z(g, S(h)_{2})\z(l, S(h)_{1}) \quad by ~(CQT2)
\\&=\z(g, S(h_{1}))\z(l, S(h_{2}))
\\&=\z(h_{1}, g)\z(h_{2}, l). \quad by ~(\ref{e6.5})
\end{align*}
We prove now (CQT3).
\begin{align*}
h_{1}g_{1}&\z(h_{2}, g_{2})
\\&=h_{1}g_{1}\v(h_{2}\rh g_{2})\quad by ~(\ref{e6.4})
\\&=(id\o \v\o \v)(h_{1}g_{1}\o h_{2}\rh g_{2}\o 1)
\\&=(id\o \v\o \v)(h_{1}(g\o 1)_{[-1]}\o h_{2}\tr (g\o 1)_{[0]})
\\&=(id\o \v\o \v)((h_{1}\tr (g\o 1))_{[-1]}h_{2}\o (h_{1}\tr (g\o 1))_{[0]})\quad by ~(\ref{y1})
\\&=(id\o \v\o \v)((h_{1}\rh g\o 1)_{[-1]}h_{2}\o (h_{1}\rh g\o 1)_{[0]})
\\&=(id\o \v\o \v)((h_{1}\rh g)_{1}h_{2}\o (h_{1}\rh g)_{2}\o 1)
\\&=(h_{1}\rh g)h_{2}
\\&=\z(h_{1}, g_{1})g_{2}h_{2}.\quad by ~(\ref{e6.7})
\end{align*}
It is easy to check that (CQT4) and (CQT5) hold.

The converse is Theorem $\ref{T6.1}$ and Proposition $\ref{P6.2}$. This completes the proof.
\end{proof}
As a corollary we have:
\begin{corollary}
Let $H$ be a Hopf algebra with a bijective antipode. Then, for $H_{4}\in \mathcal{LR}(H)$, the braiding $\psi_{H_{4}, H_{4}}$ is a symmetry if and only if $H$ is commutative.
\end{corollary}
\begin{proof}
If the braiding satisfies $\psi_{H_{4}, H_{4}}^{2}=id$, then by $(\ref{e6.4})$ $\z(k, g)=\v(k\rh g)=(\v\o \v)(k\tr (g\o 1))=(\v\o \v)(k_{1}gS(k_{2})\o 1)=\v(g)\v(k)$ for any $k,g\in H$. Thus by Theorem $\ref{T6.3}$ $(H, \v\o \v)$ is a cotriangular Hopf algebra, which by (CQT3) implies that $H$ is commutative.

Conversely, assume that $H$ is commutative, for any $k,l, g, h\in H$, we have
\begin{align*}
\psi_{H_{4}, H_{4}}&(k\o l\o g\o h)
\\&=(k\o l)_{[-1]}\tr (g\o h)_{(0)}\o (k\o l)_{[0]}\tl (g\o h)_{(1)}
\\&=k_{1}\tr (g\o h_{1})\o (k_{2}\o l)\tl h_{2}
\\&=k_{1}gS(k_{2})\o h_{1}\o k_{3}\o  h_{2}lS(h_{3})
\\&=k_{1}S(k_{2})g\o h_{1}\o k_{3}\o  lh_{2}S(h_{3})\quad by ~ H~is ~commutative
\\&=g\o h \o k\o l.
\end{align*}
It is clear that the braiding $\psi_{H_{4}, H_{4}}$ is a symmetry.
\end{proof}
If we consider $H\o \Bbbk$, by Theorem \ref{T6.3}, we generalize the another important result in \cite{CW98}.
\begin{corollary}
Let $H$ be a Hopf algebra with a bijective antipode, and assume that $(H, \rh, \D)\in \!^{H}_{H}\mathcal{YD}$, where $\D$ is usual comultiplication and $\rh$ is a left action on $H$. Then $\psi_{H, H}$ is a symmetry if and only if there exists a braiding $\z:H\o H\ra \Bbbk$ so that $(H, \z)$ is cotriangular Hopf algebra. And then $\z(k, g)=(\v\o\v)\psi(k\o  g)$, for any $k, g \in H$. That is,
\begin{align*}
k\rh g=\z(k, g_{1})g_{2}.
\end{align*}
\end{corollary}

\section*{Acknowledgements}

The second author thanks the financial support of the National Natural Science Foundation of China (Grant No. 11871144)  and the NNSF of Jiangsu Province (No. BK20171348).

\section*{ Date Availability Statement}

Data sharing is not applicable to this article as no new data were created or analyzed in this study.

\end{document}